\documentclass{article}[12pt]
\usepackage{amsmath,amsfonts,amsthm,amssymb,amsbsy,amstext,amscd,amsxtra,multicol}
\usepackage{cancel}
\usepackage{algorithm}
\usepackage[noend]{algpseudocode}
\usepackage{verbatim}
\usepackage{tikz}
\usetikzlibrary{automata,positioning}
\usepackage{multicol}
\usepackage{graphicx}
\usepackage[colorlinks,urlcolor=blue]{hyperref}
\usepackage[stable]{footmisc}
\usepackage{ dsfont }
\usepackage{wrapfig}
\usepackage{xparse}
\usepackage{ifthen}
\usepackage{bm}
\usepackage{color}
\usepackage{subcaption}
\usepackage{graphicx}
\usepackage{indentfirst}
\usepackage{enumitem}
\usepackage{xcolor}
\definecolor{xkchen}{rgb}{1.0, 0.75, 0.0}

\usepackage{xcolor}
\usepackage{hyperref}
\definecolor{linkcolor}{HTML}{799B03} 
\definecolor{urlcolor}{HTML}{799B03} 

\usepackage{amsmath,amsfonts,amsthm,amssymb,amsbsy,amstext,amscd,amsxtra,multicol}
\usepackage{cancel}
\usepackage{algorithm}
\usepackage[noend]{algpseudocode}
\usepackage{verbatim}
\usepackage{tikz}
\usetikzlibrary{automata,positioning}
\usepackage{multicol}
\usepackage{graphicx}
\usepackage[colorlinks,urlcolor=blue]{hyperref}
\usepackage[stable]{footmisc}
\usepackage{ dsfont }
\usepackage{wrapfig}
\usepackage{xparse}
\usepackage{ifthen}
\usepackage{bm}
\usepackage{color}
\usepackage{subcaption}
\usepackage{graphicx}
\usepackage{indentfirst}
\usepackage{enumitem}
\usepackage{xcolor}
\definecolor{xkchen}{rgb}{1.0, 0.75, 0.0}

\usepackage{xcolor}
\usepackage{hyperref}
\definecolor{linkcolor}{HTML}{799B03} 
\definecolor{urlcolor}{HTML}{799B03} 
\usepackage[left=2cm,right=2cm,top=2cm,bottom=2cm]{geometry}

\def \X {\mathbf{X}}
\def \Y {\mathbf{Y}}

\def \W {\widetilde{W}}

\def \range {\mathbf{range}}

\newtheorem{theorem}{Theorem}

\newtheorem{lemma}[theorem]{Lemma}
\newtheorem{assumption}[theorem]{Assumption}
\newtheorem{definition}[theorem]{Definition}
\newtheorem{remark}[theorem]{Remark} 
\makeatletter
\def\BState{\State\hskip-\ALG@thistlm}
\makeatother

\usepackage[utf8]{inputenc} 
\usepackage[T1]{fontenc}    
\usepackage{hyperref}       
\usepackage{url}            
\usepackage{booktabs}       
\usepackage{amsfonts}       
\usepackage{nicefrac}       
\usepackage{microtype}      
\usepackage{xcolor}         

\title{Adaptive Stepsize Selection in    Decentralized   Convex Optimization}

\author{%
  Ilya Kuruzov\\
  Innopolis University\\  
  \texttt{kuruzov.ia@phystech.edu}.  \\
  \and Xiaokai Chen\\\
  Purdue University\\\
  \texttt{chen4373@purdue.edu}. \\
  \and Gesualdo Scutari\\\
  Purdue University\\\
  \texttt{gscutari@purdue.edu}. \\
    \and
    Alexander Gasnikov\\
  Innopolis University
  \\\
  \texttt{gasnikov@yandex.ru} \bigskip \\
}

 \date{May 10, 2025}
\usepackage{setspace}
\begin{document}

\maketitle 

\begin{abstract}  
We study decentralized optimization where multiple  agents  minimize the average of their (strongly) convex, smooth losses over a communication graph. Convergence of the existing decentralized methods  generally  hinges on an apriori, proper   selection of the  stepsize.   Choosing this value is notoriously delicate: (i)  it demands global knowledge from all the agents of the  graph's connectivity and every local smoothness/strong-convexity constants--information they rarely have;  (ii) even with perfect information,   the worst-case tuning forces an overly small stepsize, slowing convergence in practice;   and (iii) large-scale trial-and-error tuning is prohibitive. This work   introduces a decentralized algorithm that is fully adaptive in the choice of the agents' stepsizes, without any global information and using only neighbor-to-neighbor communications--agents need not even know whether the problem is strongly convex.  The algorithm retains strong guarantees: it converges at  \emph{linear} rate when the losses are strongly convex and at \emph{sublinear} rate otherwise, matching the best-known rates of (nonadaptive) parameter-dependent methods.  
 \end{abstract}

\section{Introduction}  
Consider the following multiagent optimization problem: 
\begin{equation} 
\label{eq:init_problem}
    \min_{x \in \mathbb{R}^d} \,  \frac{1}{m}\sum\limits_{i=1}^m f_i(x),\tag{P} 
\end{equation}
where  each $f_i:\mathbb{R}^d\to \mathbb{R}$ is the loss function accessible only to   agent $i$, assumed to be      (strongly) convex and smooth (i.e., with   Lipschitz continuous gradient).  Agents communicate through a fixed, undirected, and connected mesh network, represented as a graph $\mathcal{G}$, and interactions occur exclusively among immediate neighbors without a central coordinating server. 

Problem~\eqref{eq:init_problem} arises in numerous fields including signal processing, multiagent systems, communications, and  decentralized machine learning, where data are  distributed across multiple locations. 

The literature includes various decentralized methods for solving \eqref{eq:init_problem} on mesh networks; we refer readers to the overview papers   \cite{Nedic_Olshevsky_Rabbat2018,HOngSPMagazine,9084356,Xin_Pu_Nedic_Khan_2020} and book \cite{Sayed-book}  for comprehensive reviews. Typically, existing methods rely on conservative stepsize choices, whose values depend    on global parameters such as the  Lipschitz constants of agents' functions and  network spectral properties. Since this global information is inaccessible to the agents, manual tuning is often employed, leading to unpredictable, conservative,  and problem-dependent performance.   This suggests the following  question:

 \begin{center}
{\it Can one design decentralized algorithms that adaptively select agents' stepsizes using only neighbor-to-neighbor communications, without any knowledge of the problem structure (not even  whether the local functions are convex or strongly convex) or the network topology?}
\end{center}

 Despite recent attempts to introduce adaptivity into decentralized methods \cite{Nazari_Tarzanagh_Michailidis_2022, chen2023convergence, li2024problem, Zhou24, our_nips_paper, Chen_CDC24, Notarstefano24, Ghaderyan_Werner_2025}, fully addressing this challenge remains an open problem, as documented next.

 \subsection{Related works}\label{sec:related-works}

\textbf{Adaptive centralized methods:} Recent years have seen increased interest in adaptive  optimization algorithms for centralized and federated setups. These methods include established techniques such as line-search approaches  \cite{nocedal2006numerical}, Polyak's stepsize \cite{polyak1969minimization}, Barzilai-Borwein's   stepsize \cite{BarzilaiBorwein1988}, as well as  recent advancements that   estimate the local curvature of the cost function \cite{Malitsky2019AdaptiveGD,Malitsky_Mishchenko_2024,Latafat_23b,Zhou24}.   Additional adaptive gradient methods tailored for machine learning include AdaGrad \cite{duchi2011adaptive}, Adam \cite{kingma2014adam}, AMSGrad \cite{Reddi2018OnTC}, NSGD-M \cite{pmlr-v119-cutkosky20b}, and their variants \cite{Orabona19,Ward20}. Some methods have been adapted to federated architectures \cite{50448,Li2022OnDA,Chen_Li_Li_2020}, but remain unsuitable for mesh networks due to reliance on central servers.  

{\bf Adaptive decentralized   methods:} The landscape of adaptive {\it decentralized} methods remains limited \cite{Nazari_Tarzanagh_Michailidis_2022, chen2023convergence, li2024problem, Zhou24, our_nips_paper, Chen_CDC24, Notarstefano24, Ghaderyan_Werner_2025}.  Approaches like \cite{Nazari_Tarzanagh_Michailidis_2022,chen2023convergence,li2024problem} address stochastic, online, and smooth problems, typically employing gradient normalization. These methods generally assume global Lipschitz continuity--hence, they are not applicable to  strongly convex instances of \eqref{eq:init_problem})--simplifying parameter-free convergence via standard stepsizes like $\mathcal{O}(1/\sqrt{k})$, and often still rely on some problem-specific parameter knowledge. The authors in \cite{Notarstefano24} proposes a Port-Hamiltonian framework ensuring parameter-free convergence in {\it centralized} settings. Convergence (global asymptotic stability) of decentralized implementations is ensured    for specific graph structures (e.g., cycle or fully connected graphs) or under {\it graph-dependent} stepsize constraints, calling for local knowledge of some network-dependent parameters. No explicit convergence rate expression is provided.

  The works  \cite{Emiola22, HU2021104151, Gao:BB} adapt the Barzilai-Borwein (BB) stepsize to gradient tracking algorithms \cite{sun2019distributed, di2016next, nedich2016achieving}. These methods have limited convergence guarantees: (i) agents' losses must be quadratic, as the BB rule may otherwise fail to converge {\it even in the centralized setting}; (ii) the BB-generated stepsizes must remain {\it uniformly bounded}, a condition that  generally is not satisfied or hard to check; (iii)  approaches by \cite{Gao:BB, HU2021104151} necessitate multiple communication rounds per iteration and global knowledge of network and optimization parameters, making them non-adaptive.  More recently, \cite{Chen_CDC24,Ghaderyan_Werner_2025} developed  an adaptive stepsize procedure for   gradient-tracking algorithms, inspired by the curvature estimation approach in \cite{Malitsky2019AdaptiveGD} for centralized methods. 
The method in \cite{Chen_CDC24} converges only to a neighborhood of the optimum and is numerically unstable on sparse networks  \cite{Ghaderyan_Werner_2025}. These limitations were addressed in \cite{Ghaderyan_Werner_2025}, which ensures exact convergence provided the adaptive stepsize remains below a given  threshold depending on network and optimization parameters.  Since this threshold is not known to the agents, the method is not practically implementable in decentralized settings where global information is unavailable. Additionally, both methods  \cite{Chen_CDC24,Ghaderyan_Werner_2025} lacks a convergence rate analysis, leaving their  theoretical competitiveness uncertain. 

      \begin{figure}  
      \begin{center}
    \includegraphics[scale=0.25]{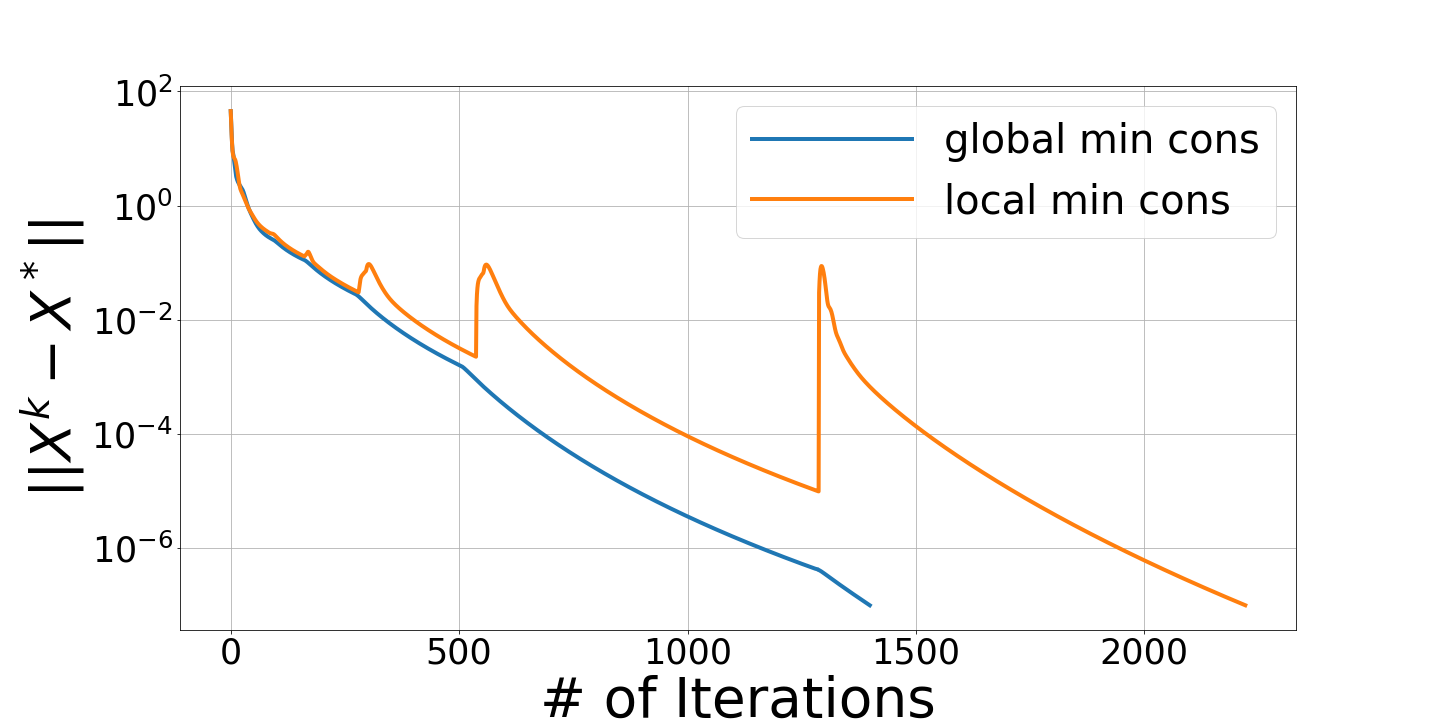}
      \end{center} 
  \caption{  Algorithm from \cite{our_nips_paper} using local and global min-consensus, applied to  ridge regression.}
  \label{fig:local_global_min_consensus}
  \end{figure}
  Recently, \cite{our_nips_paper} proposed an adaptive decentralized method for strongly convex instances of \eqref{eq:init_problem}, achieving linear convergence and demonstrating improved theoretical and practical performance over existing non-adaptive algorithms. This method employs a backtracking procedure at each agent to adaptively select stepsizes {\it locally}, followed by a \textit{global} min-consensus step synchronizing agents' local stepsizes. This global consensus is implemented using flooding protocols based on low-power wide-area network (LPWAN) technologies such as LoRa (low power, long-range) \cite{Kim_Lim_Kim_2016, Janssen_BniLam_Aernouts_Berkvens_Weyn_2020}. Alternatively, a \textit{local} min-consensus rule is proposed, which exchanges information only with {\it immediate} neighbors. While practically more feasible, this local approach results in \textit{weaker} theoretical guarantees, requiring stronger assumptions like bounded iterates, and   yields   \textit{non-monotonic} convergence trajectories. Consequently, decentralized termination criteria become hard to implement, and  error trajectories frequently show transient fluctuations or ``spikes''--see Fig.~\ref{fig:local_global_min_consensus}.\vspace{-.1cm}

 \subsection{Major contributions} 
The literature review above shows that existing adaptive decentralized algorithms  cannot tackle the question raised in this paper. We fill this gap. Our major contributions are as follows.\smallskip 
 
\textbf{Algorithm design:}
We propose a {\it fully} adaptive decentralized algorithm   that operates using only neighbor-to-neighbor communication. Each agent tunes its primal and dual stepsizes via a local backtracking rule and   then synchronizes them with its immediate neighbors through a lightweight local-minimum consensus that passes only {\it scalar} (hence negligible) messages. The procedure requires {\it zero} prior knowledge of optimization parameters (e.g., smoothness or   strong-convexity constants)    or network characteristics (degree, diameter, spectral gap). Agents need  not even   know whether their  cost function is merely convex or strongly convex. Consequently, the algorithm runs out-of-the-box--without any manual tuning or user intervention. 

On the technical side, our method departs substantially from \cite{our_nips_paper}. That algorithm adapts stepsizes with a {\it primal-only} backtracking rule followed by a {\it global} (network-wide) min-consensus, yet it lacks any mechanism to tame the dual-space imbalances that heterogeneous primal stepsizes create. We address this shortcoming by introducing explicit {\it dual} stepsize variables and a {\it joint primal–dual tracking scheme} that continuously re-aligns the two sets of stepsizes, compensating in real time for imbalances produced by independently chosen primal stepsizes, yet preserving adaptivity. This architecture eliminates the costly flooding step required in \cite{our_nips_paper} while preserving its convergence guarantees using only neighbor-to-neighbor communication.
A further dividend of our design  is a novel {\it running-diameter tracker} that estimates the ``effective'' graph diameter on the fly. Updated   via a scalar local-min consensus, it delivers the synchronizing power of a diameter-based global consensus practically at  a fraction of the communication overhead.

\smallskip
\textbf{Convergence analysis:} For   strongly convex  $f_i$'s, 
the algorithm converges linearly, matching the best nonadaptive rates. Crucially, the smoothness and strong-convexity constants   appearing in the bound are evaluated only over the {\it convex hull} of the iterates--a direct benefit of the  stepsize adaptivity--so the dependence on problem parameters is far milder than in worst-case global bounds in the literature. For merely convex losses, we   guarantee sublinear convergence. To our knowledge, this is the first adaptive decentralized method with proven guarantees in the purely convex setting. 

 \subsection{Notation} Throughout the paper  we use the following notation. 
 Capital letters denote matrices. Bold capital letters represent matrices where each row is an agent's variable, e.g., $\mathbf{X} = [x_1, \ldots, x_m]^\top$. For such matrices, the $i$-th row is denoted by the corresponding lowercase letter with the subscript $i$;  e.g., for $\X$, we write   $x_i$ (as column vector).   Let $\mathbb{S}^m$, $\mathbb{S}^m_+$, and $\mathbb{S}^m_{++}$ be the set of $m\times m$ (real) symmetric, symmetric positive semidefinite, and symmetric positive definite matrices, respectively; $A^\dagger$ denotes the Moore-Penrose pseudoinverse of $A$. The eigenvalues of   $W\in \mathbb{S}^m$ are ordered in nonincreasing order, and denoted by   $\lambda_1(W)\geq   \dots \geq \lambda_m(W)$.   
    We denote:    $[m]=\{1,\ldots,m\}$, for any integer $m\geq 1$; $[x]_+\!\!:=\max(x,0)$, $x\in \mathbb{R}$;  $1_m\in \mathbb{R}^m$  is the  vector of all ones;  $I_m$ (resp. $0_m$) is the $m\times m$ identity (resp. the $m\times m$ zero) matrix; the   information on the dimension is omitted when not necessary; $\texttt{null}(A)$ (resp. $\texttt{span}(A)$) is the nullspace (resp. range space) of the matrix $A$. Let $\langle X,Y\rangle:=\texttt{tr}(X^\top Y)$,  for any $X$ and $Y$ of suitable size (\texttt{tr}$(\bullet)$) is the trace operator;  and  $\|X\|_M:= \sqrt{\langle M X,X\rangle}$, for any symmetric, positive definite $M$  and $X$  of suitable dimensions. We still use  $\|X\|_M$ when $M$ is positive semidefinite and $X\in \texttt{span}(M)$.
 \section{Setup and Preliminaries} 
 We study Problem~\eqref{eq:init_problem}  under the following assumptions.   
\begin{assumption}
\label{assump:smooth}
{\bf (i)} Each function  $f_i$ in \eqref{eq:init_problem} is   $L$-smooth and $\mu$-strong convex on $\mathbb{R}^d$, for some $L\in (0,\infty)$ and  $\mu\in [0,\infty)$;   
 {\bf (ii)} each agent $i$ has access only to its own function $f_i$, but   $\mu$ and $L$; and {\bf (iii)} when $\mu=0$, Problem \eqref{eq:init_problem} has a nonempty solution set. 
 \vspace{-0.1cm}
\end{assumption}
Notice that   when $\mu=0,$  we simply require convexity.  We remark that for simplicity we assumed global smoothness (and strong convexity), however our results extend to functions that are only locally smooth (and strongly convex).
\begin{assumption}\label{ass-graph}
   Agents are embedded in a communication network      modeled
as an undirected, static, connected graph  $\mathcal{G}=\left([m], \mathcal{E}\right)$, where  $\mathcal{V}=[m]$ is the
set of  agents  and 
$(i,j)\in \mathcal{E}$ if there is   
 communication link (edge) between   $i$ and $j$. Let $\mathcal N_i:=\{j:\,|\, (i,j)\in \mathcal E,\,\text{for some } i\in [m]\}\cup \{i\}$ denote the set of immediate neighbors of agent $i$ (including agent $i$ itself). 
\end{assumption}

  \subsection{  The decentralized adaptive method \cite{our_nips_paper}}\label{sec:our_neurips_alg}
For strongly convex instances of  \eqref{eq:init_problem} ($\mu>0$), \cite{our_nips_paper} introduced the first {\it adaptive} decentralized algorithm that achieves  provable convergence without requiring   agents to know  optimization or network parameters. Each agents keeps a local replica $x_i\in \mathbb{R}^d$ of the global  variable $x\in \mathbb{R}^d$ along with  dual variables $y_i\in \mathbb{R}^d$. Let us stack the primal  and dual variables  into   $\mathbf{X} := [x_1, \ldots, x_m]^\top \in \mathbb{R}^{m\times d}$ and $\mathbf{Y} := [y_1, \ldots, y_m]^\top \in \mathbb{R}^{m\times d}$;  let   $\nabla F(\X):=[\nabla f_1(x_1), \cdots, \nabla f_m(x_m)]^\top$; and define the diagonal matrices of the agents' stepsizes  $\Theta^k = \mathrm{diag}(\theta_1^k, \ldots, \theta_m^k)$. The algorithm   updates    $\X^k$, $\Y^k$ and   $\theta_i^k$  as follows: starting from $(\X^0, \Y^0)\in \mathbb{R}^{m\times d}$ and $\Theta^0>0$, for any $k\geq 0,$  
\begin{equation}\label{eq:FBS_III}    \begin{aligned}
     \X^{k+1/2} &=W\, \X^k,\quad
       \Y^{k+1/2} = W \,\left(\Y^k+ \nabla F(\X^{k+1/2})\right),\smallskip\\
        \X^{k+1}& =\X^{k+1/2}- \Theta^{k}\cdot \Y^{k+1/2},\\
        \Y^{k+1} &= \Y^{k+1/2} + (I-W)\left(\Theta^{k}\right)^{-1}\X^k- \nabla F(\X^{k+1/2}).
    \end{aligned}
\end{equation} 
The stepsizes   
     are updated via   backtracking  followed by one of the two min-consensus schemes:
\begin{equation}\label{eq:stepsize}
    \theta_{i}^{k} =\left\{\begin{array}{ll}
       \displaystyle{\min_{j\in [m]}}\,  {\overline{\theta}_{j}^{k}}  & (\texttt{global min-consensus});\medskip   \\
  \displaystyle{ \min_{j\in \mathcal{N}_i} }\,  {\overline{\theta}_{j}^{k}}    & (\texttt{local min-consensus});
 \end{array}\right. \end{equation}
with 
	 $$ {\overline{\theta}_{j}^{k}} =\texttt{Backtracking}\left(\theta^{{k-1}}, f_i, x^{k+1/2}_i, -y^{k+1/2}_i, \gamma^{k-1}, \delta\right),$$
and the backtracking procedure  detailed below, wherein $\delta>0$ and $\gamma>1$ are given constants.   
\begin{algorithm}[h!]
\caption*{\texttt{Backtracking}\,($\theta$, $f$, $x$, $y$, $\gamma$, $\delta$)}
		\begin{algorithmic}
			\State $\theta^+ := \gamma \theta$\,;
              $x^+:=x+\theta^+ \,y\,$; and set $t=1$;
			\While{$f(x^+) > f(x)  + \left\langle \nabla f (x), x^+  - x \right\rangle + \frac{\delta}{2\theta^+} \|x^+-x \|^2$} 
			\State $\theta^+ \leftarrow  ({1}/{2})\theta^+$ and    $x^+ := x+\theta^+  y;$  set  $t\leftarrow t+1$;
			\EndWhile		
			\Return $\theta^+$.
\end{algorithmic}\label{alg:backtracking} 
\end{algorithm}

In  \eqref{eq:FBS_III}, $W\in  \mathbb{R}^{m\times m}$   represents the  gossip matrix, satisfying  the following  standard property  in the decentralized optimization literature \cite{Nedic_Olshevsky_Rabbat2018,Sayed-book,9084356}: $W\in \mathcal{W}_{\mathcal{G}}$, where $\mathcal{W}_{\mathcal{G}}$ is defined as follows.  
\begin{definition}[Gossip matrices] Let $\mathcal{W}_{\mathcal{G}}$ denote   the set of     matrices $\W=[\W_{ij}]_{i,j=1}^m$ that satisfy the following properties:  {\bf (i) } (compliance with $\mathcal{G}$) $\W_{ij}>0$ if $(i,j)\in \mathcal{E}$; otherwise $\W_{ij}=0$. Furthermore, $\W_{ii}>0$, for all $i\in [m]$;   and {\bf (ii)} (doubly stochastic) $\W\in \mathbb{S}^m$ and $\W 1_m=1_m.$\smallskip 
\end{definition}

Algorithm (\ref{eq:FBS_III}) operates as follows. Each agent:  (i)  performs two rounds of communication to compute intermediate variables $\X^{k+1/2}$ and $\Y^{k+1/2}$ via gossiping; (ii)   applies local backtracking on its function $f_i$ at $x_i^{k+1/2}$ along $-y_i^{k+1/2}$ to compute $\overline{\theta}_i^{k}$; and (iii) 
  executes a min-consensus protocol to establish the final stepsize $\theta_i^{k+1}$. 
In the \textit{global} min-consensus mode, agents execute a {\it network-wide} min-operation and adopt the  
common stepsize value  $\min_{j \in [m]} \overline{\theta}_j^{k}$.  
In contrast, the  \textit{local} min-consensus rule  $\min_{j \in \mathcal{N}_i} \overline{\theta}_j^k$ exchanges information only with immediate neighbors, but resulting in much weaker convergence guarantees--see discussion in Sec.~\ref{sec:related-works} and Fig.~\ref{fig:local_global_min_consensus}. 

The next section addresses these shortcomings, devising an adaptive scheme relying solely on one-hop communications yet preserving convergence guarantees of Algorithm~\eqref{eq:FBS_III} with global min-consensus.


 \vspace{-0.2cm}

\section{A New Fully-Decentralized Design} 
\label{sec:alg_design}
  
Using~(\ref{eq:FBS_III}) equipped with   the global min-consensus   as a benchmark, our first step is to understand the reasons behind the performance gap between its local and global min-consensus variants.

  We compare the one-step updates of Algorithm~(\ref{eq:FBS_III}) from  an     arbitrary  iterate $(\X^k, \Y^k)$ with $\Y^k\!\! \in\!\! \texttt{range}(I-W)$,  $(\X^k, \Y^k)$, using two distinct step-size strategies. The first one, denoted $(\X^{k+1}_{\texttt{lc}}, \Y^{k+1}_{\texttt{lc}})$, uses \textit{heterogeneous} step sizes  $\Theta^{k}= \text{diag}(\theta^{k}_1\dots \theta^k_m)>0$, with   $\theta_{\min}^{k}(:=\min_{i\in [m]} \theta_i^{k})\neq  \theta^{k}_{\max}(:=\max_{i\in [m]} \theta_i^{k})$. The second instance, denoted $(\X^{k+1}_{\texttt{gb}}, \Y^{k+1}_{\texttt{gb}})$,  uses  {\it uniform}   step sizes  
$\Theta^{{k}}=\theta^{{k}}_{\min} I$. 
The constructions   $(\X^{k+1}_{\texttt{lc}}, \Y^{k+1}_{\texttt{lc}})$ and $(\X^{k+1}_{\texttt{gb}}, \Y^{k+1}_{\texttt{gb}})$  capture   precisely the updates in  Algorithm~(\ref{eq:FBS_III}) using the    local-and  global-min consensus strategies (\ref{eq:stepsize}), respectively.  Let $(\X^\star,\Y^\star)$ denote an arbitrary fixed point of the algorithm, using global min-consensus;  it must be $\X^\star=1 (x^\star)^\top$, with $x^\star$ being an optimal solution of \eqref{eq:init_problem} \cite{our_nips_paper}, and $1\in \mathbb{R}^m$ is the vector of all ones. 

  \begin{lemma}\label{lm_iterates-error-original} 
In the setting above and Assumptions \ref{assump:smooth}-\ref{ass-graph}, the following holds:  \begin{equation}\label{eq:error_same_pi} \|\X^{k+1}_{\texttt{gb}} - {\X}^{k+1}_{\texttt{lc}}\|+\|\Y^{k+1}_{\texttt{gb}} - {\Y}^{k+1}_{\texttt{lc}}\|= \mathcal{O} \left(\left(\theta_{\max}^{k}-\theta_{\min}^{k{}}\right)\max\left(\|\X^\star\|, R^k\right)\right),\end{equation}  where    $\textcolor{black}{R^k:=} \textcolor{black}{\|\X^k-\X^\star\|+\|\Y^k-\Y^\star\|}$. 
Furthermore,   if   $\|\X^\star\|\geq 2m^{1.5} R^k $,  then  
\begin{equation}\label{eq:error_same_pi_lower_bound} \|\Y^{k+1}_{\texttt{gb}} - {\Y}^{k+1}_{\texttt{lc}}\|= \Omega \left(\left(\theta_{\max}^{k{}}-\theta_{\min}^{k}\right)\max\left(\|\X^\star\|, R^k\right)\right). \end{equation}  
  \end{lemma}
Eq.~\eqref{eq:error_same_pi_lower_bound}   indicates that whenever $\theta_{\min}^{k{}} \neq \theta_{\max}^{k{}}$--the typical scenario when employing  the local min-consensus--the dual trajectory ${\Y}^{k+1}_{\texttt{lc}}$  {can} deviate from ${\Y}^{k+1}_{\texttt{{gb}}}$, even if $R^k=0$, i.e., when the algorithm has reached the optimal solution at iteration $k$. Enforcing $\theta_{\min}^k = \theta_{\max}^k$ for every $k$   using only neighbor-to-neighbor communications  is generally infeasible and would, in any case, significantly reduce the algorithm's adaptivity.
This issue arises because of the presence of the term $\|\X^\star\|$ in the $\max$ expression in \eqref{eq:error_same_pi}-\eqref{eq:error_same_pi_lower_bound}, which persists even if $R^k$ approaches zero (or is sufficiently small). 

To overcome this limitation and  preserve    stepsize adaptivity, we propose using {\it separate}  stepsizes in the dual variable updates   in   \eqref{eq:FBS_III}, denoted by  $\Pi^k = \mathrm{diag}(\pi_1^k, \dots, \pi_m^k)$. The modification of \eqref{eq:FBS_III} reads:  
\begin{equation}\label{eq:FBS_Pi}    \begin{aligned}
     \X^{k+1/2} &=W\, \X^k,\quad
       \Y^{k+1/2} = W \,\left(\Y^k+ \nabla F(\X^{k+1/2})\right),\smallskip\\
        \X^{k+1}& =\X^{k+1/2}- \Theta^{k}\cdot \Y^{k+1/2},\\
        \Y^{k+1} &= \Y^{k+1/2} + (I-W)\left(\Pi^{k}\right)^{-1}\X^k- \nabla F(\X^{k+1/2}).
    \end{aligned} 
\end{equation}

 We proceed  designing  neighbor-only primal/dual stepsize updates. 
 
 We begin deriving the counterpart of  Lemma~\ref{lm_iterates-error-original}.  Starting from   $(\X^k, \Y^k)$,    compare one step of \eqref{eq:FBS_Pi} with {\it heterogeneous} stepsizes $\Theta^{k}$ and $\Pi^{k}$ to Algorithm~\eqref{eq:FBS_III} using the uniform stepsize $\Theta^{k{}}=\theta_{\min}^{k} I$, yielding iterates $(\X^{k+1}_{\pi},\Y^{k+1}_{\pi})$ and $(\X^{k+1}_{\texttt{gb}},\Y^{k+1}_{\texttt{gb}})$, respectively. We have the following, where $(\X^\star,\Y^\star)$ denotes an arbitrary fixed point of Algorithm (\ref{eq:FBS_III}) using uniform stepsizes.  
  

\begin{lemma}\label{lm_iterates-error-new}
In the setting above and Assumptions \ref{assump:smooth}-\ref{ass-graph}, the following holds:  
\begin{equation}\label{eq:Delta_error_new}\Delta^{k+1}= \mathcal{O} \left(\Big(\big|\theta^{k}_{\min}-\pi^{k}_{\min}\big|+ \big(\theta^{k}_{\max}-\theta^{k}_{\min}\big)\Big) R^k+ \left(\pi^{k}_{\max}-\pi^{k}_{\min}\right)\max\left(\|\X^\star\|, R^k\right)\right),\end{equation}
where $\Delta^{k+1}:=\|\X^{k+1}_{\texttt{gb}} - {\X}^{k+1}_{{\pi}}\|+\|\Y^{k+1}_{\texttt{gb}} - {\Y}^{k+1}_{{\pi}}\|$, $\pi_{\min}^k:=\min_i\pi_i^k$, and $\pi_{\max}^k:=\max_i\pi_i^k$.
\end{lemma}
 Eq.~\eqref{eq:Delta_error_new} highlights that the discrepancy $\Delta^{k+1}$ is now controllable even when $\theta_{\min}^{k} \neq \theta_{\max}^{k}$, provided $\pi_{\min}^{k} \approx \pi_{\max}^{k}$.  This suggests the following  desideratum for the   updates  of  $\Theta^k$ and $\Pi^k$: 
 \begin{itemize}
     \item[\textbf{(c1)}]  Maintain $\Pi^k \approx \pi^k I$ for almost all iterations, and   some suitably adjusted scalar $\pi^k$; and
     \item[\textbf{(c2)}] Ensure that $|\pi^k_{\min}-\theta_{\min}^k|$ and $\theta^k_{\max}-\theta_{\min}^k$ remain  uniformly bounded, possibly sufficiently small across all iterations  $k$.
 \end{itemize}    
     \textbf{(c1)} guarantees that the second term in \eqref{eq:Delta_error_new} remains negligible  while \textbf{(c2)} promotes adherence of iterates~\eqref{eq:FBS_Pi}   to those of the benchmark  \eqref{eq:FBS_III} that uses  global min-consensus. 

Based on the above, we propose the following  updating procedures for the primal and dual stepsizes. \smallskip 

  \textbf{Primal stepsize updates:} Given $\theta_i^k$, each agent $i$  tunes ${\overline{\theta}_i^k}$ via the backtracking procedure (Sec.~\ref{sec:our_neurips_alg}), then runs the {\it local} min-consensus step (\ref{eq:stepsize}) on    ${\overline{\theta}_i^k}$, producing $\theta_i^{k}$. 
This method adapts stepsizes to local curvatures of the losses $f_i$, still  controlling the discrepancy $\theta^k_{\max} - \theta^k_{\min}$. 
\smallskip 

  \textbf{Dual stepsize updates:}   Ideally,  $\Pi^k = \theta_{\min}^k I$, for all $k$,  would  fulfill  \textbf{(c1)} and \textbf{(c2)}, but $\theta_{\min}^k$ is {\it unknown} locally. 
To address this, we design a neighbor‑only protocol  tracking $\theta_{\min}^k$. We hinge on the following   properties of the 
backtracking  search, which hold  irrespectively of the algorithm trajectory.
\begin{lemma} \label{lemma:step_size_increase_decrease} Given an $L$-smooth function $f$,   arbitrary  
 sequences $\{x^k, y^k\}_k \subseteq \mathbb{R}^d \times \mathbb{R}^d$ and $\{\gamma^k\}_k\subset \mathbb{R}_+$, and $\theta^0>0$ and $\delta>0$, let   $\theta^{k+1} = \texttt{Backtracking} (\theta^k, f, x^k, y^k, \gamma^{k}, \delta)$, $k\in \mathbb{N}_0$. 
 The following hold, for every $k$:  
  (a) either  $\theta^{k+1}<\theta^{k}$ or   $\theta^{k+1}=\gamma^k\theta^k$;  and (b) the number of decreases of stepsize value, $\mathcal{I}_{\theta}(k) =\{j\leq k\,|\,\theta^{j+1}<\theta^j\}$, is bounded by   $|\mathcal{I}_{\theta}(k)|=\mathcal{O}(\sum_{j=1}^k \log\gamma^j).$ \vspace{-0.2cm}
\end{lemma}

Initially, assume each agent knows the graph diameter $d_{\mathcal{G}}$ and receives $\theta_{\min}^k$ precisely every $d_{\mathcal{G}}$ iterations, that is,  at   $k\equiv0\mod d_{\mathcal{G}}$. The   horizon of $d_{\mathcal G}$
  is because the min‑consensus protocol  
  converges in at most  $d_{\mathcal G}$
  steps.  
 In the intermediate iterations  \textcolor{black}{$k^\prime=k+1,\ldots, k+d_{\mathcal G}-1$}, 
 if 
 $\{\theta_{\min}^{k^\prime}\}$ does not decrease,  by Lemma~\ref{lemma:step_size_increase_decrease}(a),  each agent can compute   $\theta^{\textcolor{black}{k^\prime}}_{\min}$ directly without   communicating, by 
$\theta^{k^\prime}_{\min}=\gamma^{k^\prime-1} \theta^{k^\prime-1}_{\min} $.  
  This suggests   the following   update for  $\Pi^k$: \begin{equation}\label{eq:Pi_update_ideal}{\Pi^k=\theta^k_{\min}I,\,\,\,  \forall k\equiv 0\text{ mod } d_{\mathcal G},\quad \text{and}\quad 
      \Pi^{k^\prime}= \gamma^{k^\prime-1 }\Pi^{k^\prime-1},\,\,\, \forall k^\prime={k+1},\ldots, k+d_{G}-1.}
  \end{equation}

This ensures \textbf{(c1)} always, and \textbf{(c2)} except for the iterations where $\theta_{\min}^{\color{black} k^\prime}$ decreases; 
Lemma~\ref{lemma:step_size_increase_decrease}(b) shows that such drops are rare--at most 
 $\mathcal{O}(\sum_{j=1}^k \log\gamma^j)$ within $\{0,\ldots, k\}$, for any  $k\in \mathbb{N}_0$.  
  We will  choose    $\{\gamma^k\}$ exactly to   minimize such occurrences, guided by the convergence analysis (Sec.~3). 

The rule \eqref{eq:Pi_update_ideal}      subsumes that every agent   obtains  $\theta_{\min}^k$, which  would require a {\it global} min‑consensus protocol. To avoid that,   we propose 
    tracking   $\theta_{\min}^k$   introducing  {\it local surrogates} collected in $\widetilde{\Theta}^k=\text{diag}(\widetilde{\theta}_1^k,\dots,\widetilde{\theta}_m^k)$. Each agent $i$ updates  $\widetilde{\theta}_i^k$ performing a {\it local} min-consensus step:   $\forall k\in \mathbb{N}_0$, let
\begin{equation}\label{eq:local_min_d}\widetilde{\theta}_i^{k}=\min_{j\in \mathcal N_i} \theta_j^{k},\,\,\text{if}\,\, k \equiv 1 \text{ mod }d_{\mathcal G};\quad \text{and} \quad  \widetilde{\theta}_i^{k^\prime}=\min_{j\in \mathcal N_i} \gamma^{k^\prime{-1}}\widetilde{\theta}_j^{k^\prime-1},  \,\,\text{if}\,\,   k^\prime={k+1},\ldots, k+d_{G}-1.\end{equation}
Every $d_{\mathcal G}$ iterations, at $k\equiv 1\,\text{mod } d_{\mathcal G}$, the min-consensus step \eqref{eq:local_min_d} is ``reset'' with the current primal stepsize value ${\theta}_i^k$, for each agent $i$. The above procedure has the following properties.   \textbf{(i)} \textit{Consensus}: after at most   $d_{\mathcal G}$ hops, a consensus is reached,   $\widetilde{\theta}^{k+d_{\mathcal{G}}{-1}}_i=\widetilde{\theta}^{k+d_{\mathcal{G}}{-1}}_j$, for all $i,j\in [m]$, $i\neq j$,   and $k\in \mathbb{N}_0$. \textbf{(ii)} \textit{Recovery of  $\theta_{\min}^k$:}   If  $\{\theta_{\min}^{k\prime}\}$ does not decrease over $k^\prime=k+1, \ldots, k+d_{\mathcal G}-1$, it must be 
$\widetilde{\theta}^{k+d_{\mathcal{G}}{-1}}_i={\theta}_{\min}^{k+d_{\mathcal{G}}{-1}}$, for all $i\in [m]$--hence agents recover   ${\theta}_{\min}^{k+d_{\mathcal{G}}{-1}}$ {\it locally}. If  $\{\theta_{\min}^{k\prime}\}$ drops, $\widetilde{\theta}^{k}_i$ will be still consensual at $k+d_{\mathcal{G}}{-1}$,  approximating  ${\theta}_{\min}^{k+d_{\mathcal{G}}-1}$.

 Using \eqref{eq:local_min_d}, we replace \eqref{eq:Pi_update_ideal} with the following:  \begin{equation}\label{eq:Pi_update_tilde_lambda}{\Pi^k=\widetilde{\Theta}^k,\quad k\equiv 0\text{ mod } d_{\mathcal G}\quad \text{and}\quad 
      \Pi^{k^\prime}= \gamma^{k^\prime-1 }\Pi^{k^\prime-1},\quad \forall k^\prime={k+1},\ldots, k+d_{G}-1.}
  \end{equation}
Since $\widetilde{\Theta}^{k}$   is consensual at $k+d_{\mathcal G}\!-\!1$, the   re-initializations of $\Pi^k$ at $k\equiv 0\text{ mod } d_{\mathcal G}$ in \eqref{eq:Pi_update_tilde_lambda} yield   $\Pi^k=\pi^k I$, that is  {\it consensual}   dual stepsizes. If $\{\theta_{\min}^{k^\prime}\}$ does not decrease   in the previous $d_{\mathcal G}{-1}$ iterations (up to $k$), by the recovery property of \eqref{eq:local_min_d}, 
 it follows $\pi^k\!=\!\theta_{\min}^k$, matching the ideal rule \eqref{eq:Pi_update_ideal} {\it without  requiring global min-consensus}.\smallskip 


\textbf{Adaptive  estimation of the ``effective'' diameter:} 
Our  rules  \eqref{eq:local_min_d}-\eqref{eq:Pi_update_tilde_lambda} need a horizon long enough for the local‐min consensus to converge. When the {\it true} graph diameter $d_{\mathcal{G}}$ is unknown, we introduce a fully-decentralized procedure  whereby each agent $i$ keeps an estimate $d_i^k$ that replaces   $d_{\mathcal G}$ in \eqref{eq:local_min_d}-\eqref{eq:Pi_update_tilde_lambda}, and enlarges it only when strictly necessary, eventually converging to an {\it effective diameter}--the smallest integer  that still guarantees \eqref{eq:local_min_d} and \eqref{eq:Pi_update_tilde_lambda} to work exactly as intended.  Let the auxiliary variables $\widetilde{\theta}^k_i$ evolve as in \eqref{eq:local_min_d}, but with $d_{\mathcal{G}}$ replaced by the current guess $d_i^k$; the $d_i^k$'s are updated     
for all $k\in \mathbb{N}_0$ and $i\in [m]$, as 
\begin{equation}\label{eq:dim_update} d_i^{k+1}= \max_{j\in\mathcal{N}_i} \,(2 d_j^{k}),\,\,\,\text{if}\,\,\,    k\equiv 0\text{ mod } d_i^k\quad \text{and}\quad   \widetilde{\theta}^{k}_i\neq \min_{j\in\mathcal{N}_i}\widetilde{\theta}^{k}_j;\quad \text{otherwise}\quad 
      d_i^{k+1}=\max_{j\in\mathcal{N}_i} d_j^{k}.
  \end{equation}
 In words, at every iteration that is a multiple of the current horizon, $k\equiv  0 \text{ mod }d_i^{k}$, agent $i$ checks whether the local-min consensus has already converged,  $\widetilde{\theta}^{k}_i\overset{?}{=}\min_{j\in \mathcal{N}_i}\widetilde{\theta}^{k}_j$. If the test fails, the diameter estimate is doubled. This value   then feeds a one-step max-consensus  for synchronization with the other agents'.   { \bf  Why the test is sound:} 
  Take two consecutive multiples of the current horizon,
  $k^\prime$ and $k^{\prime\prime} (>k^\prime)$, so that  $k^\prime\equiv 0\text{ mod } d_i^{k^\prime}$ and $k^{\prime\prime}\equiv 0\text{ mod } d_i^{k^{\prime\prime}}$; it is     $k'' \geq k' + d_i^{k'}$. According to \eqref{eq:local_min_d},  no re-initialization of the $\widetilde{\theta}$-variables occurs in   $[k'+1, k'']$. Hence, if  the local consensus test fails at $k^{\prime\prime}$  then at least  
  $d_i^{k'}$ rounds of local min-consensus proved insufficient, and doubling is required. 
   { \bf  Guarantees:} \textbf{(i)} {\it  Finite number of failures}--A consensual value $d_i^k=d_{\mathcal G}$, for all $i$,   always passes the test; specifically, one can show that   each agent can fail it at most  $\sim \log_2 d_{\mathcal G}$. 
   \textbf{(ii)}{\it Horizon may be smaller than the true diameter}--The local–min consensus in \eqref{eq:local_min_d} can converge on a common horizon
   $d_i^k=d_j^k$, $i\neq j\in [m]$, that is  strictly smaller than $d_{\mathcal G}$. This   facilitates condition \textbf{(c2)}, as it reduces  the frequency of iterations where primal and dual step sizes differ, thus bringing the proposed fully-decentralized scheme closer to Algorithm (\ref{eq:FBS_III}) employing    global-consensus communications.  

   The overall algorithm incorporating the primal-and dual-step procedures \eqref{eq:local_min_d} and \eqref{eq:Pi_update_tilde_lambda} as well as the adaptive estimation scheme in (\ref{eq:dim_update}) of the effective diameter is summarized in  Algorithm~\ref{alg:adaptive_procedure}.

   The simple example below highlights why an adaptive horizon can outperform the conservative choice $d_{\mathcal G}$.

\begin{remark}
   Consider a line graph with 
  $m$ agents, so   $d_{\mathcal G}=m-1.$ Assume the cost functions alternate along the line:  $f_{2n+1}\equiv f$ and  $f_{2n}\equiv g$, for $n=1\dots \lfloor m/2\rfloor$; where   $f,g:\mathbb{R}^d\rightarrow \mathbb{R}$ are some convex and smooth functions.   Additionally, suppose the backtracking yields identical   stepsizes for equal functions, i.e.,   $\textcolor{blue}{\overline{\theta}}^{k}_{2n+1}\equiv \theta_f^k$ and $\textcolor{blue}{\overline{\theta}}^k_{2n}\equiv \theta_g^k$, for all $k.$  This is the case, e.g.,  for $f(x)=x^2$ and $g(x)=ax^2, a>1.$ One round of local min-consensus then is enough to synchronize all agents: $\theta^{\textcolor{blue}{k}}_i = \min(\theta^k_f, \theta^k_g),$  for all $i\in [m]$.  Hence, all   min-consensus procedures in Algorithm~\ref{alg:adaptive_procedure} converge   in one round.  This    matches the behavior of Algorithm (\ref{eq:FBS_III}) with global min-consensus! In contrast, if every agent   uses   $d_i^k=d_{\mathcal G}$, it  would introduce a  $m-1$-fold delay, potentially slowing convergence.  
\end{remark}\vspace{-0.3cm}

  \begin{algorithm}[ht!]		\caption{ }
  \noindent \textbf{Data:}  (i) Initialization:    $\X^0 \in \mathbb{R}^{m\times n}$  and  $\Y^0 = 0$;  $\theta^{-1}_i \in (0, \infty)$ (primal-stepsize),   $d_i^0:=d^0$ for some $d^0\in\mathbb{N}$ (local diameter estimate);   $ \pi^{-1}_i, \widetilde{\theta}^{-1}_i \in (0, \infty)$  (dual stepsizes and auxiliary variables).  (ii) Backtracking   parameters: $\delta \in (0,1]$, and    $\{\gamma^k\}_k\subseteq [1,\infty)$. (iv) Gossip matrix  $W:= (1-c)I_m+c \widetilde{W}$, with $ \widetilde{W} \in \mathcal{W}_{\mathcal{G}}$, and $c\in (0,1/2]$.  Set the iteration index  $k = 0$.		\smallskip 
		\begin{algorithmic}
           
           \State \texttt{(S.1) Communication step:} Agents updates primal and dual variables via gossiping:  
            $$  \X^{k+1/2}  =W\, \X^k \quad\text{and}\quad
       \Y^{k+1/2} = W \,\left(\Y^k+ \nabla F(\X^{k+1/2})\right);$$
			\State \texttt{(S.2) Update of the primal stepsizes}: Each agent $i\in [m]$ updates their primal step sizes using a local backtracking procedure followed by a min-consensus step:  
   $$\overline{\theta}_i^{\textcolor{black}{k}} ={\texttt{Backtracking}}\left(\theta_i^{\textcolor{black}{k-1}}, f_i, x^{k+1/2}_i, y^{k+1/2}_i, \gamma^{k-1}, \delta\right);\quad \theta^{\textcolor{black}{k}}_i = \min_{j\in \mathcal{N}_i} \overline{\theta}^k_j.$$
			\State \texttt{(S.3) Update of the  dual stepsizes and auxiliary variables}: Each agent $i\in [m]$ updates their axillary variables and dual stepsizes:
           $$\ \widetilde{\theta}^{{k}}_i=\begin{cases}\min_{j\in \mathcal N_i} \theta_j^{k},&\text{if $\textcolor{black}{k\equiv 1\, \text{mod } d_i^{k}}$},\\ \min_{j\in \mathcal N_i} \gamma^{k}\widetilde{\theta}_j^{k},&\text{else,}\end{cases}\quad \pi^{\textcolor{black}{k}}_i=\begin{cases} \widetilde{\theta}_i^{{k}},&\text{if $ {k}\equiv 0\, \text{mod } d_i^{k}$},\\ \gamma^k{\pi}^{\textcolor{black}{k-1}}_i,&\text{else.}\end{cases} $$

\State \texttt{(S.4) Graph diameter estimation:} Each agent $i$ updates its local diameter estimation:   
$${{d}_i^{k+1}}=\begin{cases}  \max_{j\in \mathcal{N}_i} (2\, d_j^{k}),&\text{if   $k\equiv  0 \,\text{mod }d_i^{k}\,$ and $\,\widetilde{\theta}^{k}_i\neq \min_{j\in\mathcal{N}_i}\widetilde{\theta}^{k}_j$}, \\ \max_{j\in \mathcal{N}_i} d_j^{k},&\text{otherwise.}\end{cases}$$

     \State \texttt{(S.5) Local updates of the primal and dual variables:} $$  
    \begin{aligned}
        \X^{k+1}& =\X^{k+1/2}- \Theta^{\textcolor{black}{k}}\cdot \Y^{k+1/2},\quad \X^{k+1/2}_\Pi = W\,(\Pi^{\textcolor{black}{k}})^{-1} \X^k,\\
        \Y^{k+1} &= \Y^{k+1/2} +  (\Pi^{\textcolor{black}{k}})^{-1} \X^k-\X^{k+1/2}_\Pi-\nabla F(\X^{k+1/2}).\\
    \end{aligned}$$
			\State \texttt{(S.6)}    If a termination criterion is not met,  $k\leftarrow k+1$ and go to step \texttt{(S.1)}.
\end{algorithmic}\label{alg:adaptive_procedure}
\end{algorithm}

 \vspace{-.1cm}

\section{Convergence Analysis}\label{sec:convergence}\vspace{-.2cm}
Before diving into the convergence analysis, we introduce some preliminary definitions and facts. 

\begin{lemma}\label{lemma:optimality_conditions}
    Consider Problem~\eqref{eq:init_problem} under Assumption~\ref{assump:smooth}; let  $(\X^{\star},\Y^{\star})$ be a fixed point of Algorithm~\ref{alg:adaptive_procedure}; and let $F:\mathbb{R}^{m\times d}\rightarrow \mathbb{R}$ be the augmented function defined as $$F(\X):=\frac{1}{m}\sum_{i=1}^m f_i(x_i),\quad x_i\in \mathbb{R}^d,\,\,i\in [m].$$  Then, the following hold: \begin{itemize}
        \item[\bf (i)] $\X^{\star}=\mathbf{1}(x^{\star})^\top,$ where $x^{\star}$ is a solution of~\eqref{eq:init_problem}; \item[\bf (ii)] $\Y^{\star}\in\range(I-W)$ and $\nabla F(\X^{\star})+\Y^{\star}=0.$ 
    \end{itemize}     
\end{lemma}
Every fixed point of Algorithm~\ref{alg:adaptive_procedure} yields consensual rows of $\X^\star$, being a solution of \eqref{eq:init_problem}. 

\subsection{The case of strongly convex $f_i$}\label{sec:convergence_text_scvx} We introduce the following  merit function evaluated along the iterates $\{\X^k,\Y^k\}$ of Algorithm~\ref{alg:adaptive_procedure}:  
$$V^{k}:=\|\X^k-\X^{\star}\|^2+(\theta^{k-1}_{\min})^2\|\Y^k-\Y^{\star}\|^2_M,$$
where $(\X^\star, \Y^\star)$ is a fixed point of the algorithm (with $\X^\star= 1(x^\star)^\top$ and $x^\star$ being the unique solution of \eqref{eq:init_problem}); $M:=c^{-1}(I-\W)^{\dagger}-I$, with $(I-\W)^{\dagger}$ denoting  the Moore-Penrose pseudoinverse of $I-\W$;       $\|\cdot\|_M$ is the norm induced by the matrix $M$;  and $\theta_{\min}^k=\min_{i\in [m]}\theta_i^k$. Notice that $V^k$ is a valid merit function: $V^k=0$ if and only if $\X^k=\X^\star$ and $\Y^k=\Y^\star$ (due to $\Y^k\in \range{(I-W)}=\range(I-\widetilde{W})$). 
 
Linear convergence of Algorithm~\ref{alg:adaptive_procedure} applied to  strongly convex instances of \eqref{eq:init_problem} is summarized next. 
 \begin{theorem}
\label{theorem:strong_convex} Given Problem~\eqref{eq:init_problem} under Assumption~\ref{assump:smooth}, with $\mu>0$ and unique   solution $x^\star\in \mathbb{R}^d$;   let $\{(\X^k, \Y^k)\}$ be the iterates generated by Algorithm~\ref{alg:adaptive_procedure}. Additionally,   choose  $\{\gamma^k\}$ as 
    $\gamma^k\leq \big((k+\beta_1)/(k+1)\big)^{\beta_2}$, for all $k$ and some   $\beta_1\geq 1, \beta_2>0.$ Then 
 {$V^k\leq \varepsilon$}, 
for all $k\geq N_\varepsilon$, with 
\begin{equation} \label{eq:N_case2} N_\varepsilon ={\mathcal O}\left(\frac{1}{\delta}\frac{\kappa{\log_2d_{\mathcal{G}}}}{(1-c(1-\lambda_m(\W)))^2 c(1-\lambda_2({\W}))} (\log(\max\{\|\X^{\star}\|^2,V^0\}/{\varepsilon}) + d_{\mathcal{G}})\right).\end{equation}
Here  $\kappa$ is the condition number of each $f_i$  restricted to  the convex hull of   $\{x^\star, \{x_i^k, x^{k+1/2}_i\}_{k=0}^{N_\varepsilon}\}$, and $\mathcal{O}$ hides the dependence on $\beta_1$ and $\beta_2$. 
\end{theorem}

  Because  $1-c(1-\lambda_m(\W))>1-2c$, the linear rate (\ref{eq:N_case2}) yields   (omitting poly-log factors)  $$N_\varepsilon=\widetilde{\mathcal{O}}\left(\frac{\kappa}{1-\lambda_2(\W)}\log \Big(\frac{1}{\varepsilon}\Big)+d_{\mathcal G}\right).$$    Quite remarkably, this matches the complexity of the adaptive method in \cite{our_nips_paper}, which requires \emph{global} min-consensus, whereas Algorithm~\ref{alg:adaptive_procedure} uses only one-hop exchanges; the extra term  $d_\mathcal G$ is somehow expected, due to the min-consensus procedures.   
  The complexity bound above is also on par with non-adaptive gradient-tracking schemes such as~\cite{nedich2016achieving}, SONATA~\cite{sun2019distributed}, and \cite{LiNA-linear18}. Yet our condition number $\kappa$ is \emph{local}--evaluated on the convex hull of the trajectory--typically much smaller than the {\it global} condition number that governs the performance of those methods.  

  The sequence $\{\gamma^k\}_{k=1}^\infty$  (each $\gamma^k\geq 1$) used  in Line 1 of the backtracking procedure,    encourages \emph{non-monotone} growth, enabling larger stepsizes between consecutive line-searches. Any choice with $\gamma\downarrow 1$ and $\prod_{k=1}^\infty \gamma^k=\infty$ is suitable. In Theorem~\ref{theorem:strong_convex}  we suggested a practical option we found effective.   For maximal simplicity,   one may instead set  $\gamma^k=1$, for all $k$, foregoing this extra parameter.

\subsection{The case of (nonstrongly) convex $f_i$}\label{sec:convergence_cvx}




In this setting, we will use the following merit function to monitoring progresses of the algorithm towards optimality:  for any given fixed point $(\X^\star,\Y^\star)$ of Algorithm~\ref{alg:adaptive_procedure},   let   $$\mathcal{M}(\X) := \max\left(\delta \|\X\|_{I-W}^2, F(\X)-F(\X^{\star})+\langle\Y^{\star}, \X\rangle\right),$$
where $\delta>0$ is as defined in the algorithm, and $\X^\star=1(x^\star)^\top$.  Note that    $\mathcal{M}(\X)=0$ if and only if   $\X=\mathbf{1} x^\top$, for some $x$ such that $(1/m)\sum\limits_{i=1}^mf_i(x)\leq (1/m)\sum\limits_{i=1}^mf_i(x^*)$ (such an $\X$ is such that $\langle\Y^{\star}, \X\rangle=0$, due to $\Y^\star\in \range \,(I-W)$--see Lemma~\ref{lemma:optimality_conditions}). 

Convergence of  Algorithm~\ref{alg:adaptive_procedure}  is summarized next, where $(\X^\star,\Y^\star)$ is any given fixed point.
 
 \begin{theorem}
\label{theorem:singular_convex} Given Problem~\eqref{eq:init_problem} under Assumption~\ref{assump:smooth}, with $\mu=0$, let $\{(\X^k, \Y^k)\}$ be the iterates generated by Algorithm~\ref{alg:adaptive_procedure}. Choose  $\{\gamma^k\}$ as 
    $\gamma^k\leq \big((k+\beta_1)/(k+1)\big)^{\beta_2}$, for all $k$ and some   $\beta_1\geq 1, \beta_2>0$. Further, let $R>0$ and $c_\theta$>0 such that    
    $R\geq \|\X^{k+1}-\X^\star\|$,   $R^2\geq \max(V^0, \|\X^{\star}\|^2)$; and   $c_{\theta}\geq L\theta^k_{\min}$. 
    Then, $$\mathcal{M}(\widehat{\X}^k) \leq \varepsilon,\quad \text{with } \quad  \widehat{\X}^{k}:=\frac{1}{k}\sum\limits_{t=1}^{k}\X^t,$$ for all $k\geq N_\varepsilon$, and 
$$N_{\varepsilon}=\mathcal{O}\left(\frac{c_{\theta}d_{\mathcal{G}}\tilde{L}  R^2 \log1/\varepsilon }{\varepsilon} + \frac{d_{\mathcal{G}} \log d_{\mathcal G}\tilde{L}R^2}{\varepsilon}\right),$$
where $\tilde{L}$ is the Lipschitz constant of each function $f_i$  restricted to  the convex hull of   $\{x^\star, \{x_i^k, x^{k+1/2}_i\}_{k=0}^{\infty}\}$,   and $\mathcal{O}$ hides the dependence on $\beta_1$ and $\beta_2$.   
\end{theorem}
Theorem~\ref{theorem:singular_convex} guarantees a complexity of $\mathcal{O}(1/\varepsilon)$, aligning with the complexity of nonaccelerated, nonadaptive (first-order) decentralized optimization methods, e.g., \cite{li2017decentralized, shi2015extra} applied to \eqref{eq:init_problem} ($\mu=0$). Notably, the proposed method eliminates the need for parameter tuning. Furthermore, similar to the strongly convex scenario, the complexity depends on parameters defined over the convex hull formed by the algorithm’s trajectory. To  our knowledge, this is the first {\it adaptive} decentralized method for (nonstrongly) convex problems with provable convergence. 
\smallskip 

 \textbf{On the compactness of $\{\X^k,\Y^k\}$:}  Notice that Theorem~\ref{theorem:singular_convex} requires compactness of the iterates $\{\X^k,\Y^k\}$. This can be ensured by   suitable choices of the sequence $\{\gamma^k\}$--two options are the following.
\begin{itemize}
    \item[\textbf{(i)}]    Choose   $\gamma^k=1$, for all $k\geq K$ and some finite $K\in \mathbb{N}_+$.  In such a case, 
    one can show that there exists some finite iteration index, after which  all primal (and thus dual) stepsizes  do  not increase, resulting in constant stepsizes over time and uniform across the agents. 
    \item[\textbf{(ii)}]  Given  $\widetilde{R}\in\mathbb{R}_{++}$,   let us  introduce additional local binary variables $h_i^k\in\{0,1\}$ (initialized $h_i^0=1$), defined as  
    \begin{equation}\label{eq:h_update}h_i^k = \begin{cases}0, &\quad \text{if }  \max\left(\|x^k_i-x^0_i\|, \theta^{k-1}_i\|y^k_i-y^0_i\|\right)\geq \widetilde{R},\\
       \displaystyle{\min_{j\in\mathcal{N}_i}} h_i^{k-1}, &\quad \text{otherwise.}\end{cases}\end{equation}
        This extra   variables monitor    the condition  
        $\max\left(\|x^k_i-x^0_i\|, \theta^{k-1}_i\|y^k_i-y^0_i\|\right)\leq \widetilde{R}.$
        
        Equipped with these variables,  Algorithm~\ref{alg:adaptive_procedure} changes as follows: \begin{itemize}
            \item[\bf (a)] Add   step \texttt{(S.0)} (before \texttt{(S.1)}) wherein (\ref{eq:h_update}) is performed;  and 
               \item[\bf (b)] Define
        $$\widetilde{\gamma}^{k-1}_i = 1 + h_i^{k-1} ({\gamma}^{k-1}-1),$$
     and call the    \texttt{Backtracking} procedure  (step \texttt{(S.2)}), replacing     $\gamma^{k-1}$ with   $\widetilde{\gamma}^{k-1}_i$.  
        \end{itemize}
        This will guarantee that all the iterates generated by this variant of the algorithm are bounded. 
        
        

\end{itemize}

 \section{Numerical Results } \label{section:num_exp}  
In this section, we compare the proposed method against the adaptive decentralized algorithm from \cite{our_nips_paper} (considering both local and global min-consensus implementations) and the well-known EXTRA algorithm \cite{shi2015extra}, whose stepsize is fine tuned for the best practical performance on the problems and network simulated. EXTRA  serves as a benchmark representing decentralized algorithms that rely on  stepsizes selected with full knowledge of both network topology and optimization parameters.

\subsection{Strongly convex quadratic} As strongly convex instance of    \eqref{eq:init_problem}, we consider a quadratic optimization problem,  with each $f_i(x)=\|A_i x - b_i\|^2$. Elements of $A_i\in\mathbb{R}^{h\times n}$ and $b_i\in\mathbb{R}^h$, with $h=110$ and $n=100$, are  generated as i.i.d. from the standard normal distribution.  
We simulated three different network graphs, of size   $m=20$, namely: (i) line graph; (ii) Erd\H{o}s--R\'{e}nyi graph with edge probability $p=0.1$ (modeling a sparse network); and (iii) Erd\H{o}s--R\'{e}nyi graph with edge probability $p=0.5$ (modeling a fairly  connected network).  In our algorithm, we choose $\gamma^k=(k+2)/(k+1)$.

\begin{figure}[h]  \hspace{-2.2cm}
\includegraphics[scale=0.14]{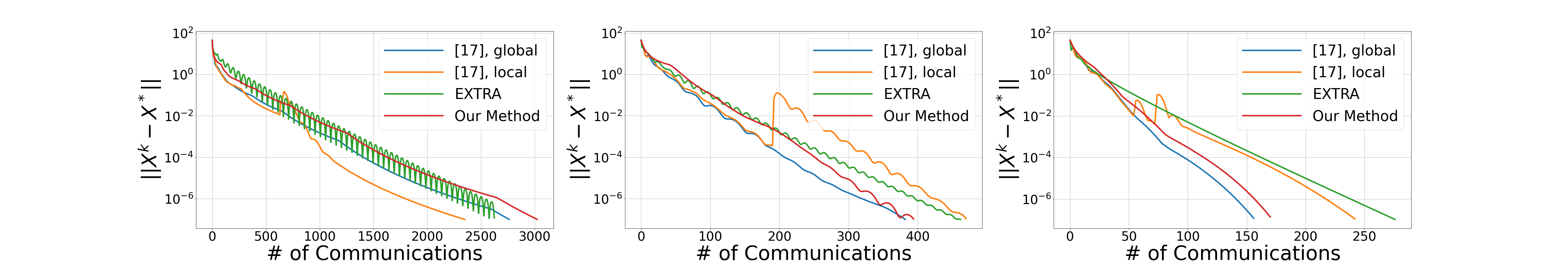}
\caption{\textbf{Strongly convex quadratic} program  on different graphs: (left) Line graph; (middle) Erdős-Rényi Graph with edge activation probability $p=0.1$; (right) Erdős-Rényi Graph with  $p=0.5$}\label{fig:ridge_compar} 
\end{figure}

The comparison among the algorithms described above is illustrated in Fig.~\ref{fig:ridge_compar}, where we plot the distance of the iterates   from the optimal solution versus the number of  {communications}. Remarkably, the proposed method consistently compares favorably against the adaptive algorithm from \cite{our_nips_paper}  relying on \emph{global} min-consensus communications. Notably, the performance gap diminishes as network connectivity improves. Additionally, our approach demonstrates significantly greater stability compared to the local min-consensus variant from \cite{our_nips_paper} and  outperforms the EXTRA algorithm, despite requiring no prior knowledge of network topology or optimization parameters. 

\subsubsection{Dependence on the condition number}   Let us consider now $f_i(x)=\|A_ix-b_i\|^2+\lambda\|x\|^2/2$ so that changing   $\lambda$ one can control the condition number of the problem. Graphs and matrices $A_i, b_i$ are generated as described in   Figure~\ref{fig:ridge_compar},   and so is the  procedure followed to  tune  the simulated algorithms.

\begin{figure}[ht]
     \centering
     \begin{subfigure}[b]{0.32\textwidth}
         \centering
         \includegraphics[width=1.15\textwidth]{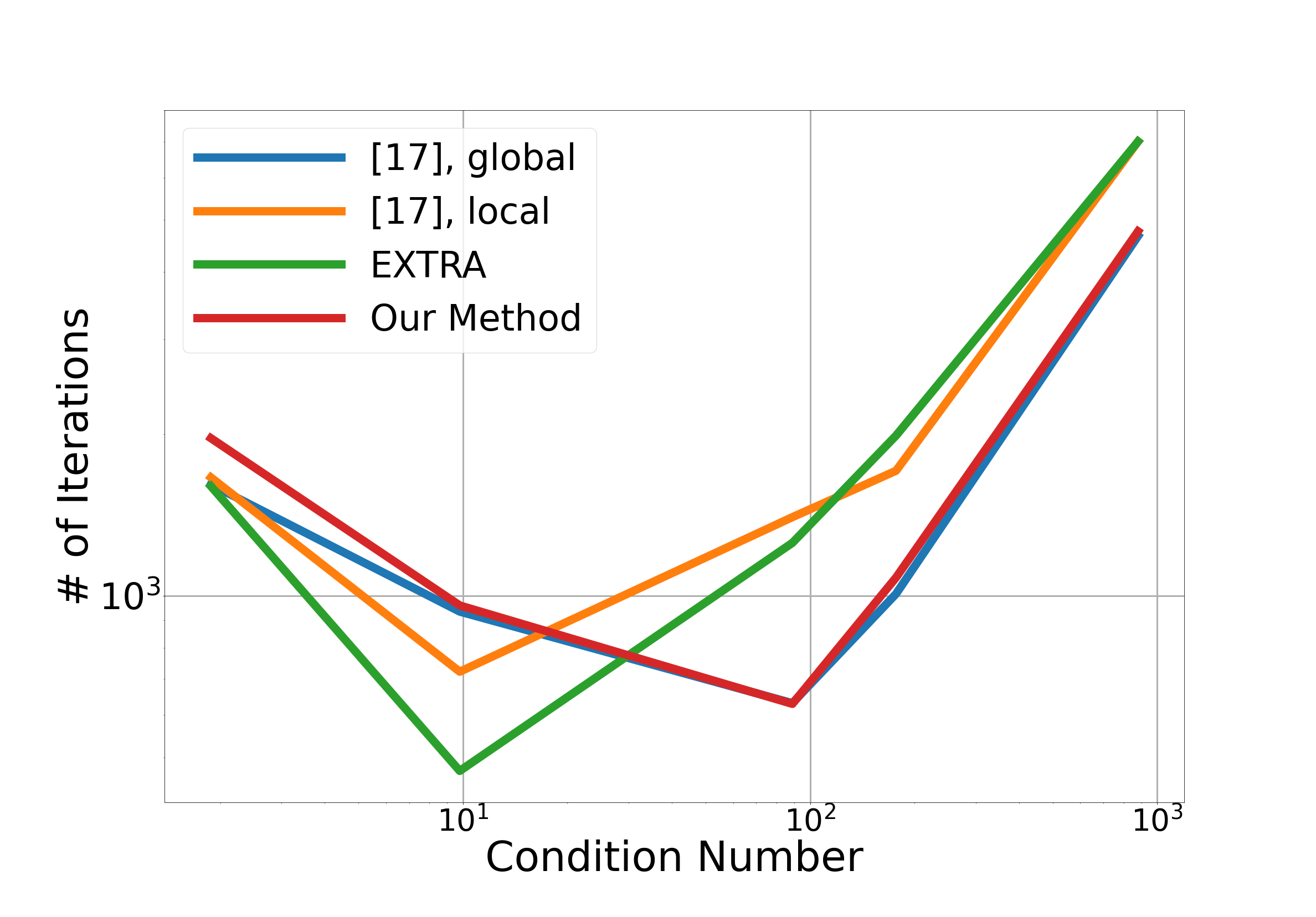}
         \caption{Line Graph}
         \label{pic1_cond_num}
     \end{subfigure}
\hfill
\begin{subfigure}[b]{0.32\textwidth}
         \centering
         \includegraphics[width=1.15\textwidth]{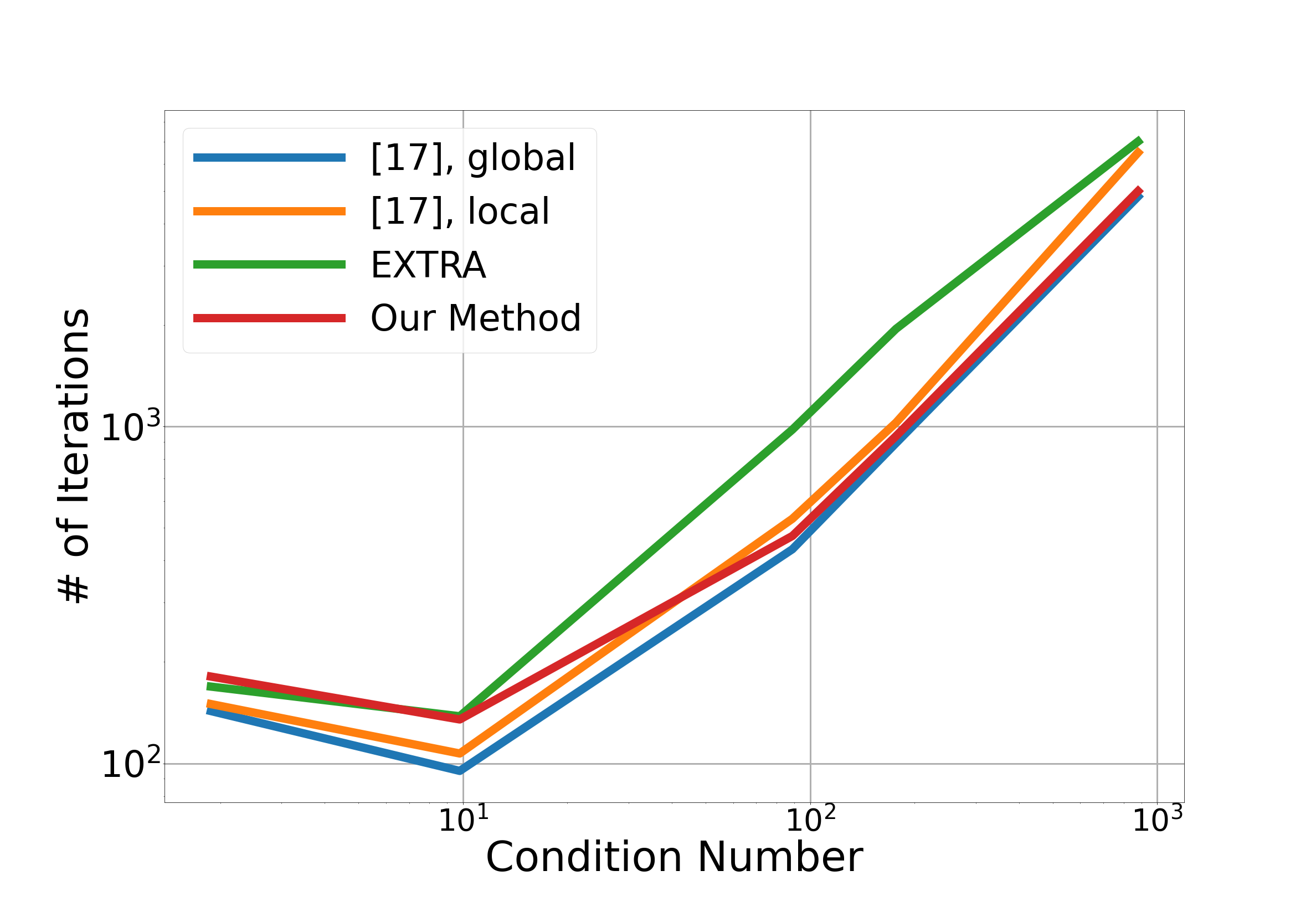}
         \caption{Erdős-Rényi Graph,  $p=0.1$}
         \label{pic2_cond_num}
     \end{subfigure}
     \hfill
     \begin{subfigure}[b]{0.32\textwidth}
         \centering
         \includegraphics[width=1.15\textwidth]{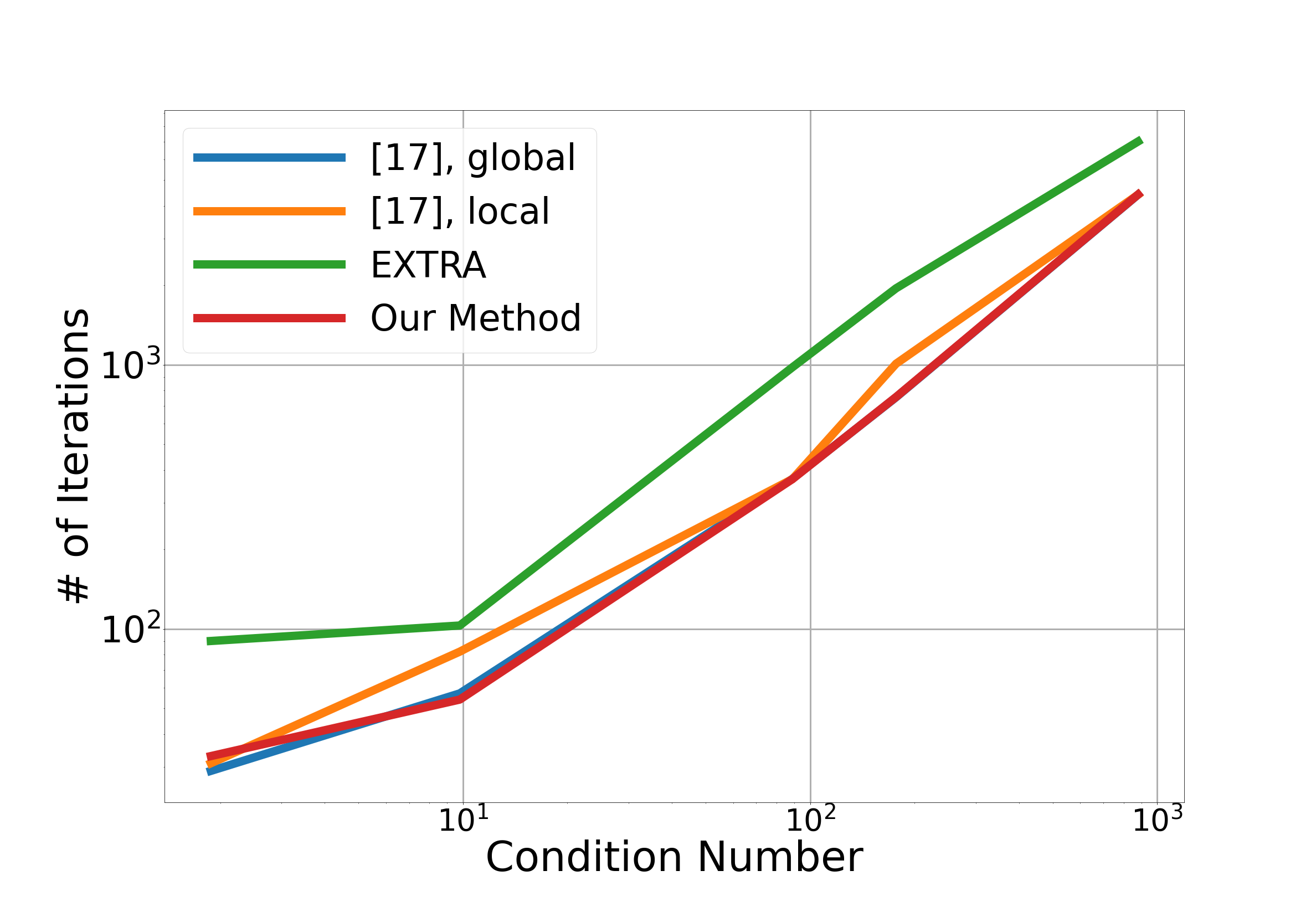}
         \caption{Erdős-Rényi Graph,   $p=0.5$}
         \label{pic3_cond_num}
     \end{subfigure}\smallskip 
     
        \caption{\small \textbf{Ridge regression}:  Number  of iterations $N$ for $\|\X^N-\X^\star\|\leq 10^{-5}$ versus the condition number  of agents' looses   on different graphs; (a) Line graph; (b) Erdős-Rényi Graph, with edge activation probability $p=0.1$; and (c) Erdős-Rényi Graph, with edge activation probability $p=0.5$.}
\label{fig:ridge_compar_cond_number} 
\end{figure}

Results are reported in Figure~\ref{fig:ridge_compar_cond_number}. Quite remarkable, the proposed  method outperforms EXTRA and \cite{our_nips_paper} using local-min consensus and performs as good as  \cite{our_nips_paper} using {\it global} min-consensus, over all graphs simulated and mostly of the condition number values.  

\subsubsection{Dependence on the graph diameter} The  algorithm complexity as stated in Theorems~\ref{theorem:strong_convex} and ~\ref{theorem:singular_convex} contains some dependence on graph diameter $d_{\mathcal G}$. It is natural then to ask how the graph diameter value affects the numerical performance of the algorithm.  We consider  the same function $f_i(x)=\|A_i x - b_i\|^2$ with  $A_i\in\mathbb{R}^{h\times n}$ with dimensions  $h=1$, $n = 100, and $ $d=100$. We simulated this problem over a line-graph with variable number of agents, hence variable diameter.   

\begin{figure}[ht]
\centering
    \includegraphics[scale=0.1]{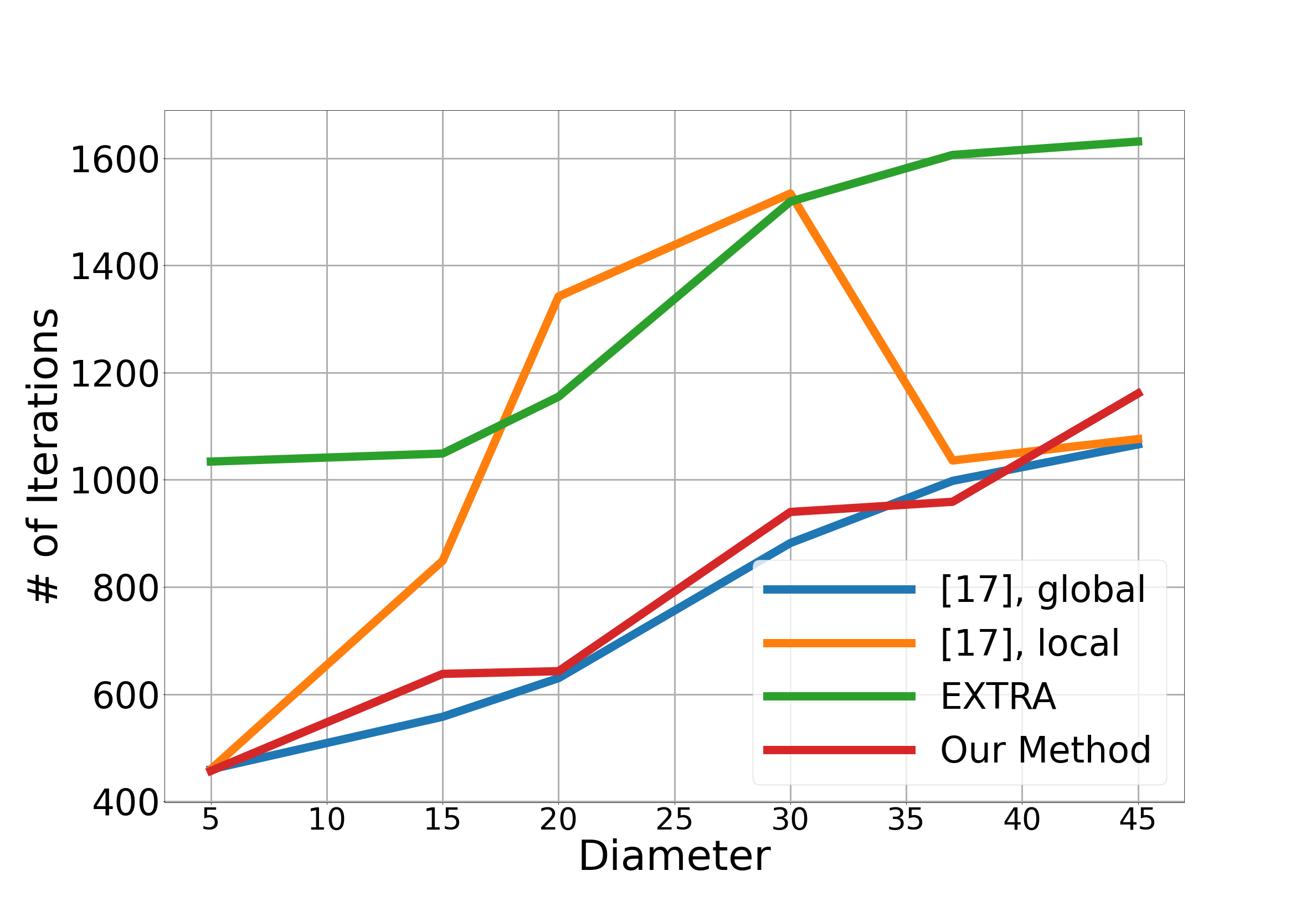}
\caption{\textbf{Strongly convex quadratic} program  on graph-line with different number of nodes: number  of iterations $N$ for $\|\X^N-\X^\star\|\leq 10^{-5}$ versus graph diameter}\label{fig:diameter} 
\end{figure}

Figure~\ref{fig:diameter} summarizes the results. Quite remarkably,  the size of the graph seems not to   affect the performance of the algorithm more than that of \cite{our_nips_paper} employing global min synchronization of the agents' primal stepsizes.  This is despite the additional dependence of the algorithm complexity on the graph diameter in Theorems~\ref{theorem:strong_convex} and ~\ref{theorem:singular_convex}.  

\subsection{Logistic regression} As (nonstrongly) convex instance of   \eqref{eq:init_problem}, we consider the decentralized logistic regression problem, corresponding to   $f_i(x)=({1}/{h})\sum_{j=1}^h\log(1+\exp(-b_{i,j}\langle f_{i,j}, a_{i,j}\rangle))$. Here,   $b_{ij}\in \{0, 1\}$, $a_{ij}\in\mathbb{R}^{200}$ are  data  problem, taken from the dataset a3a \cite{libsvm}. We distributed data across $m=20$ nodes, each owning  $h=159$ samples. We simulated the same networks as in Fig.~\ref{fig:ridge_compar}. 

\begin{figure}[h] 
	\hspace{-2.1cm}
\includegraphics[scale=0.14]{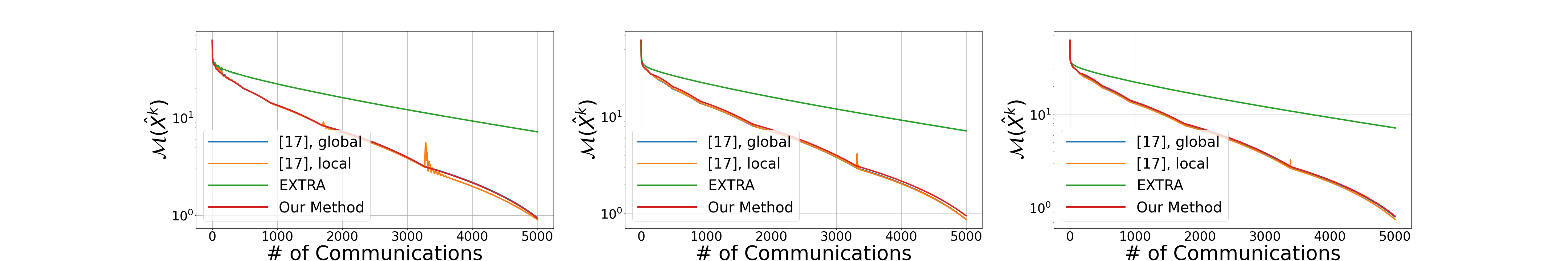}\smallskip 

\caption{\textbf{Logistic regression} on different graphs: (left) Line graph; (middle) Erdős-Rényi Graph with edge activation probability $p=0.1$; (right) Erdős-Rényi Graph with  $p=0.5$ }\label{fig:logreg_compar}
\end{figure}

Results are summarized in Fig.~\ref{fig:logreg_compar}; we plot the values of the merit function $\mathcal M$ along the ergodic iterates generated by all the algorithms versus the number of  {communications}. 
Even on (nonstrongly) convex problems, the proposed method   outperforms   EXTRA and  \cite{our_nips_paper} using  local min-consensus while being  comparable with \cite{our_nips_paper}. using  global min-consensus. In such a setting however  there is no convergence proof for \cite{our_nips_paper}. Hence the proposed method is the first {\it adaptive} decentralized algorithm for (nonstrongly) convex instances of \eqref{eq:init_problem} with provable convergence.

 \bibliography{biblio}
 \bibliographystyle{plain} 
\end{document}